\documentclass[parskip=half,bibliography=totoc]{scrartcl}

\usepackage[automark]{scrlayer-scrpage}								
\usepackage{amsfonts}												
\usepackage{amsmath}												
\usepackage{mathtools}												
\usepackage{amssymb}												
\usepackage{extarrows}												
\usepackage{dsfont} 												
\usepackage{mathrsfs}												
\usepackage{accents}												
\usepackage[T1]{fontenc}											
\usepackage[utf8]{inputenc}											
\usepackage{subcaption} 											
\usepackage[top=1in,bottom=1in,left=1.25in,right=1.25in]{geometry}	
\usepackage{footmisc}												
\usepackage{enumerate} 												
\usepackage{booktabs}												
\usepackage{authblk}												
\usepackage{abstract}												
\usepackage{soul}													
\usepackage{stmaryrd}												
\usepackage{comment}												

\usepackage{tikz-cd}												

\usepackage{amsthm}													

\usepackage[linktocpage=true,colorlinks=true,citecolor=blue,urlcolor=blue]{hyperref}				
\usepackage[noabbrev]{cleveref}										
\crefname{lem}{lemma}{lemmata}

\usepackage[noadjust]{cite}

\renewcommand{\d}{\mathrm{d}}

\newcommand{\CH}{\mathbb{C}\mathrm{H}}
\newcommand{\HH}{\mathbb{H}\mathrm{H}}

\newcommand{\C}{\mathbb{C}}
\renewcommand{\H}{\mathbb{H}}

\newcommand{\calC}{\mathcal{C}}
\newcommand{\calH}{\mathcal{H}}
\newcommand{\calL}{\mathcal{L}}
\newcommand{\calS}{\mathcal{S}}

\newcommand{\hor}{\mathrm{H}}
\newcommand{\ver}{\mathrm{V}}

\DeclareMathOperator{\Span}{span}

\DeclareMathOperator{\End}{End}

\DeclareMathOperator{\Aut}{Aut}
\DeclareMathOperator{\Isom}{Isom}

\DeclareMathOperator{\Heis}{Heis}

\DeclareMathOperator{\Rm}{Rm}

\newcommand{\h}{\mathrm{H}}

\newcommand{\norm}[1]{\lVert #1 \rVert}

\newtheoremstyle{mythm}
{}
{}
{\slshape}
{}
{\bfseries\sffamily}
{.}
{ }
{}
\newtheoremstyle{mydef}
{}
{}
{}
{}
{\bfseries\sffamily}
{.}
{ }
{}

\theoremstyle{mythm}
\newtheorem{thm}{Theorem}[section]
\newtheorem{prop}[thm]{Proposition}
\newtheorem{cor}[thm]{Corollary}
\newtheorem{lem}[thm]{Lemma}
\theoremstyle{mydef}
\newtheorem{mydef}[thm]{Definition}
\newtheorem{rem}[thm]{Remark}

\newenvironment{myproof}[1][\proofname]{
	\proof[\sffamily\upshape#1]
}{\endproof}

\clearscrheadfoot
\ihead[]{\headmark}
\ohead[]{\pagemark}
\deffootnote[1em]{0em}{1em}{%
	\textsuperscript{\thefootnotemark}%
}
\setfootnoterule{3em}

\newenvironment{numberedlist}{\begin{enumerate}[\upshape(i)]}{\end{enumerate}}

\allowdisplaybreaks 
\title{{A class of locally inhomogeneous complete quaternionic K\"ahler manifolds}}
\author{Vicente Cort\'es}
\author{Alejandro Gil-Garc\'ia}
\author{Arpan Saha}

\affil{\normalsize Fachbereich Mathematik\\
	Universit\"at Hamburg\\
	Bundesstra\ss e 55, 20146 Hamburg, Germany\\
	vicente.cortes@uni-hamburg.de, alejandro.gil.garcia@uni-hamburg.de}
\affil{\normalsize Instituto de Ciencias Matem\'aticas\\
	Nicol\'as Cabrera 13--15, Cantoblanco 28049 Madrid, Spain\\ arpan.saha@icmat.es}
\date{}

\newcommand{\HK}{\mathrm{HK}}
\newcommand{\ra}{\mathrm{a}}
\newcommand{\rb}{\mathrm{b}}

\date{\large\today}

\begin{document}

\maketitle

\begin{abstract}
We prove that the one-loop deformation of any quaternionic K\"ahler manifold in the class of c-map spaces is locally inhomogeneous. As a corollary, we obtain that the full isometry group of the one-loop deformation of any homogeneous c-map space has precisely cohomogeneity one.

\end{abstract}

\clearpage

\tableofcontents

\section{Introduction}

It is by now well known that the supergravity c-map \cite{FS} and its one-loop deformation \cite{RSV} can be used to construct a wealth of complete quaternionic K\"ahler manifolds of negative scalar curvature, see \cite{ACDM} for a geometric construction of the relevant metrics and \cite{CHM,CDS} for some completeness theorems. The local geometry of any c-map space can be encoded in a holomorphic function subject to a (non-holomorphic) non-degeneracy condition on its two-jet, the so-called holomorphic prepotential of special K\"ahler geometry. In \cite{CST21,CST22} it was shown that the one-loop deformation of one of the classical series of symmetric quaternionic K\"ahler manifolds of negative scalar curvature has precisely cohomogeneity one.

The purpose of this paper is to show that the one-loop deformation of any c-map space is locally inhomogeneous, see Theorem~\ref{main:thm}. As a consequence, we prove that the one-loop deformation of any homogeneous quaternionic K\"ahler manifold of negative scalar curvature, with exception of the quaternionic hyperbolic space (that is not a c-map space), has precisely cohomogeneity one, see Corollary~\ref{main:cor}. Our arguments rely on a general formula for the curvature tensor of any (possibly indefinite) hyper-K\"ahler manifold obtained from the rigid c-map, see Theorem~\ref{theorem:curvature_N}, as well as on special properties of the curvature tensor in the case of rigid c-map spaces associated with \emph{conical} affine special K\"ahler manifolds, see Proposition~\ref{cor:curvature_N_HZ=0}. The above information obtained for the hyper-K\"ahler manifolds in the image of the rigid c-map is transferred to crucial information about the curvature tensor of the quaternionic K\"ahler manifolds in the image of the one-loop deformed supergravity c-map in Proposition~\ref{crucial:prop}. Finally, this information is used in Proposition~\ref{main:prop} to show that the quaternionic K\"ahler manifold has a non-constant scalar-valued curvature invariant when the deformation parameter $c$ is positive. This essentially implies the claimed inhomogeneity and cohomogeneity results.

\subsection*{Acknowledgements}

We thank Danu Thung for comments and suggestions regarding an earlier draft of this paper.

This work was supported by the German Science Foundation (DFG) under Germany's Excellence Strategy  --  EXC 2121 ``Quantum Universe'' -- 390833306. 
A. S. is supported by the Spanish Ministry of Science and Innovation, through the `Severo Ochoa Programme for Centres of Excellence in R\&D' (CEX2019-000904-S).

\section{Background}

In this section, we will recall some relevant background material. In particular, we recall the definitions and properties of conical affine special K\"ahler manifolds, as well as the rigid and (deformed) local c-map constructions which assign pseudo-hyper-K\"ahler manifolds to affine special K\"ahler manifolds and a one-parameter family of quaternionic K\"ahler manifolds to conical affine special K\"ahler manifolds, respectively.

\subsection{Special K\"ahler geometry}

\begin{mydef}
	An \emph{affine special K\"ahler (ASK) manifold} $(M,g,J,\omega,\nabla)$ is a (pseudo-) K\"ahler manifold $(M,g,J,\omega )$ endowed with a flat torsion-free connection $\nabla$ such that $\nabla\omega=0$ and $\d^\nabla J = 0$.
\end{mydef}

\begin{mydef}
	A \emph{conical affine special K\"ahler (CASK) manifold} $(M,g,J,\omega,\nabla,\xi)$ is an ASK manifold $(M,g,J,\omega,\nabla)$ endowed with a complete vector field $\xi$ such that 
	\begin{itemize}
		\item $g$ is negative-definite on the span of $\xi$ and $J\xi$, and positive-definite on its orthogonal complement,
		\item $D\xi =\nabla\xi =\mathrm{id}$, where $D$ denotes the Levi-Civita associated to $g$. 
	\end{itemize}
\end{mydef}

\begin{rem}
	As $D(J\xi)= JD \xi = J$ is skew-symmetric with respect to the metric $g$, it follows that $J\xi$ is a Killing vector field. Similarly, we see that $J\xi$ is Hamiltonian.
\end{rem}

We will only consider connected CASK manifolds, so in the following all the manifolds involved will be connected. Moreover, we will always assume that the vector fields $\xi$ and $J\xi$ generate a principal $\C^\times$-action. The quotient manifold $\bar M$ does then inherit a (positive-definite) K\"ahler metric $\bar g$ and $(\bar M, \bar g)$ is called a \emph{projective special K\"ahler manifold}. The metric is obtained by K\"ahler reduction exploiting the fact that $J\xi$ is a Hamiltonian Killing vector field.

Let us define the tensor field
\begin{equation*}
	\begin{split}
		\calS&:=g^{-1}\nabla g\in\Gamma(T^*M\otimes\End(TM)).
	\end{split}
\end{equation*}
In $\nabla$-affine coordinates $q^i$, this is given by
\begin{equation}\label{eq:Christoffel}
	\calS_{ij}^k=g^{km}(\nabla g)_{ijm}=2\Gamma_{ij}^k.
\end{equation}
The following three results are well known in special K\"ahler geometry \cite{Freed,ACD2002}.

\begin{lem}\label{Ssymm_lem}
	The tensor $\calS$ on an ASK $(M,g,J,\nabla)$ manifold satisfies the following properties: \begin{enumerate}[\normalfont(a)]
		\item $\calS_XY=\calS_YX$.
		\item $g(\calS_XY,Z)=g(\calS_XZ,Y)$.
	\end{enumerate} for all vector fields $X,Y,Z\in\Gamma(TM)$.
\end{lem}

\begin{myproof}
	Since $\calS=g^{-1}\nabla g$, we have $g(\calS_XY,Z)=(\nabla_Xg)(Y,Z)$. Since $\nabla\omega=0$ and $g=\omega(\cdot,J\cdot)$, we furthermore have $g(\calS_XY,Z)=\omega(Y,(\nabla_XJ)Z)$.
	
	Part (a) follows from the fact that $(\nabla_XJ)Z-(\nabla_XJ)Z=(\d^\nabla J)(X,Z)=0$. Part (b) follows from the fact that $g$ is symmetric.
\end{myproof}

Note that in particular this means that $\nabla g$ is fully symmetric.

\begin{prop}
	Let $(M,g,J,\nabla)$ be an ASK manifold. Then we have $D-\nabla=\tfrac{1}{2}\calS$, where $D$ is the Levi-Civita connection of $g$ and $\calS=g^{-1}\nabla g$.
\end{prop}

\begin{myproof}
	Let $\widetilde{D}:=\nabla+\tfrac{1}{2}\calS$. As the Levi-Civita connection is the unique torsion-free connection preserving the metric $g$, the result will follow if we can show that $\widetilde D$ is torsion-free and metric.
	
	Since $\nabla$ is torsion-free and $\calS_XY=\calS_YX$, it follows that $\widetilde{D}$ is also torsion-free. Now let us check that $\widetilde{D}$ is metric. For $X,Y,Z\in\Gamma(TM)$ we have $$\begin{aligned}
		(\widetilde{D}_Xg)(Y,Z)&=Xg(Y,Z)-g(\widetilde{D}_XY,Z)-g(Y,\widetilde{D}_XZ)\\
		&=Xg(Y,Z)-g(\nabla_XY,Z)-g(Y,\nabla_XZ)-\tfrac{1}{2}g(\calS_XY,Z)-\tfrac{1}{2}g(Y,\calS_X Z)\\
		&=(\nabla_X g)(Y,Z)-\tfrac{1}{2}(\nabla_X g)(Y,Z)-\tfrac{1}{2}(\nabla_X g)(Z,Y)=0,
	\end{aligned}$$ where in the last step we have used that $\nabla g$ is fully symmetric.
\end{myproof}

\begin{cor}\label{cor:S_vanish}
	Let $(M,g,J,\nabla,\xi)$ be a CASK manifold. Then for all $X\in\Gamma(TM)$, we have $$\calS_\xi X=\calS_X\xi=0,\qquad\calS_{J\xi}X=\calS_XJ\xi=0.$$
\end{cor}

\begin{myproof}
	The first equation follows from $\tfrac{1}{2}\calS_X\xi=D_X\xi-\nabla_X\xi=0$ and that $\calS_XY=\calS_YX$ for all $X,Y\in\Gamma(TM)$. For the second one, we note $$\begin{aligned}
		\tfrac{1}{2}\calS_XJ\xi&=D_XJ\xi-\nabla_XJ\xi=J(D_X\xi)-J(\nabla_X\xi)-(\nabla_XJ)\xi\\
		&=-(\nabla_XJ)\xi=-(\nabla_\xi J)X=\calS_\xi JX=0,
	\end{aligned}$$ where we have used $\tfrac{1}{2}\calS=D-\nabla=-\tfrac{1}{2}J\nabla J$ in the penultimate equality.
\end{myproof}

\subsection{Rigid c-map}\label{ssec:rigid-c-map}

The rigid c-map \cite{CFG,ACD2002} assigns to each affine special (pseudo-)K\"ahler manifold $(M,g,J,\nabla)$ and, in particular, to any conical affine special K\"ahler manifold $(M,g,J,\nabla,\xi)$ of real dimension $2n+2$ a pseudo-hyper-K\"ahler manifold $(N=T^*M,g_N,I_1,I_2,I_3)$ of real dimension $4n+4$.

Using the flat torsion-free connection $\nabla$, we can identify $$TN=T(T^*M)=T^\hor N\oplus T^\ver N\cong\pi^*TM\oplus\pi^*T^*M,$$ where $\pi\colon N=T^*M\longrightarrow M$ is the canonical projection, $T^\ver N=\ker(\d\pi)$ is the vertical distribution and $T^\hor N$ is the horizontal distribution defined by $\nabla$. In particular, given a vector field $X\in \Gamma(TN)$, we will think of its horizontal component $X^\hor$ as a section of $\pi^*TM$ and its vertical component $X^\ver$ as a section of $\pi^*T^*M$. Using these identifications, we have $$g_N=\begin{pmatrix}
	g&0\\
	0&g^{-1}
\end{pmatrix},\qquad I_1=\begin{pmatrix}
	J&0\\
	0&J^*
\end{pmatrix},\qquad I_2=\begin{pmatrix}
	0&-\omega^{-1}\\
	\omega&0
\end{pmatrix},\qquad I_3=I_1I_2.$$

In the case that $M$ is CASK, the rigid c-map space $N=T^*M$ enjoys some additional properties.

\begin{prop}[{\cite[Proposition~2]{ACM}}]\label{prop:rotating}
	Let $(M,g,J,\nabla,\xi)$ be a CASK manifold and define on the associated rigid c-map space $(N=T^*M,g_N,I_1,I_2,I_3)$ the following data: $$Z:=-\widetilde{J\xi},\qquad \omega_1 :=g_NI_1,\qquad \omega_\h := \omega_1 + \d \iota_Z g_N,$$ $$f_Z := -\tfrac12 g_N(Z,Z),\qquad f_\h := \tfrac12 g_N(Z,Z),$$ where $\widetilde{J\xi}$ denotes the horizontal lift of $J\xi$ with respect to $\nabla$. Then
\begin{equation*}
	\calL_Z g_N=0,\qquad \iota_Z \omega_1 = -\d f_Z,\qquad \calL_Z \omega_2 = \omega_3,\qquad \iota_Z\omega_\h = -\d f_\h.
\end{equation*}
\end{prop}

Note that Proposition~\ref{prop:rotating} implies that $Z$ is a \emph{rotating} Killing vector field for the pseudo-hyper-K\"ahler manifold $N$. This means that it is a Killing vector field that preserves one of the K\"ahler structures $I_1$ and rotates the other two into each other.

It will also be convenient for later purposes to introduce some notation. We denote by $\H Z=\Span\{Z,I_1Z,I_2Z,I_3Z\}$ the distribution generated by the quaternionic span of $Z$, then
\begin{equation}\label{decomp:eq}
	TN=\H Z\oplus(\H Z)^\perp.
\end{equation}
The vector fields $Z,I_1Z\in\Gamma(TN)$ are horizontal and the vector fields $I_2Z,I_3Z\in\Gamma(TN)$ are vertical, with respect to the decomposition $TN=T^\hor N\oplus T^\ver N$.

\subsection{HK/QK correspondence, supergravity c-map, and twist construction}\label{ssec:sugra-c-map}

Suppose we are given a pseudo-hyper-K\"ahler manifold $(N,g_{N},I_1,I_2,I_3)$ that admits \emph{HK/QK data} i.e. a tuple $(\omega_1,Z,\omega_\h,f^c_Z,f^c_\h)$ such that
\begin{itemize}
	\item $\omega_1 :=g_NI_1$ is integral,
	\item $Z$ is a rotating Killing vector field preserving $\omega_1$ (assume for simplicity that $Z$ generates a free circle action),
	\item $\omega_\h := \omega_1 + \d \iota_{Z} g_N$,
	\item $f^c_Z$ is a nowhere vanishing function such that $\iota_{Z}\omega_1 = -\d f^c_Z$
	\item $f^c_\h := f^c_Z + g_{N}(Z,Z)$ is nowhere vanishing.
\end{itemize}

Note that there is a freedom of adding a constant to the Hamiltonian functions $f^c_Z$ and $f^c_\h$, so long as the shifted Hamiltonian functions are still nowhere vanishing. This is reflected in the superscript $c$ in $f^c_Z$ and $f^c_\h$.

Given the above data, it was shown in \cite{ACM} generalising \cite{Haydys} that we can construct a quaternionic pseudo-K\"ahler manifold $(\bar N^c,g^c_{\bar N},\mathcal Q^c)$ of non-zero scalar curvature with a circle action such that the given pseudo-hyper-K\"ahler manifold $N$ may be retrieved as a hyper-K\"ahler reduction of the Swann bundle of $\bar N^c$ by a lift of the circle action at a non-zero level set. This construction is known as the \emph{HK/QK correspondence}.

The signature of the resulting quaternionic pseudo-K\"ahler manifold and the sign of its scalar curvature depend on the signature of the pseudo-hyper-K\"ahler manifold $N$ and the signs of the functions $f_Z^c$ and $f_\h^c$. The cases when one obtains a (positive-definite) quaternionic K\"ahler metric were specified in \cite{ACM} and include the case of quaternionic K\"ahler metrics of positive scalar curvature considered by Haydys, who started with a (positive-definite) hyper-K\"ahler metric. In the following theorem we focus on the cases which yield a positive definite metric of negative scalar curvature, of relevance to the present paper.

\begin{thm}[{\cite[Corollary~2]{ACM}}]\label{tm:negative_QK}
	Let $(N,g_{N},I_1,I_2,I_3)$ be a pseudo-hyper-K\"ahler manifold of real dimension $4n+4\ge 4$ equipped with HK/QK data $(Z,\omega_1,\omega_\h,f^c_Z,f^c_\h)$ and let $P$ be a principal $S^1$-bundle such that $c_1(P)=[\omega_1]=[\omega_\h]$. Then there is a lift of the circle action on $N$ generated by $Z$ to $P\times\H^\times$, so that its quotient $\hat M$ by the lifted action carries a conical pseudo-hyper-K\"ahler structure with hyper-K\"ahler reduction $(N,g_{N},I_1,I_2,I_3)$. The conical pseudo-hyper-K\"ahler manifold $\hat M$ is the Swann bundle of a (positive-definite) quaternionic K\"ahler manifold $(\bar N^c,g^c_{\bar N},\mathcal Q^c)$ of negative scalar curvature if and only if $g_{N}$ is positive-definite and $f^c_Z>0$, or if the signature of $g_N$ is $(4n,4)$ and $f^c_Z>0$ while $f^c_\hor<0$.\footnote{Note that the functions $f^c_Z$ and $f_\h^c$ correspond to $-\frac12 f$ and $-\frac12 f_1$ in \cite{ACM}, as indicated in \cite[page~100]{CST21} (up to a typo).}
\end{thm}

Note that \emph{explicit} expressions for all of the above data, including the quaternionic K\"ahler metric, are obtained in \cite{Haydys, ACM}, \cite[Theorem~2]{ACDM}. We have however omitted these in the statement of the theorem to avoid redundancy, since we will describe the metric below using the language of Swann's twist construction.

Moreover, since there is the freedom of adding a constant term to the Hamiltonian function $f^c_Z$, Theorem~\ref{tm:negative_QK} gives us, if necessary after restricting to open sets, a one-parameter family of quaternionic K\"ahler manifolds of fixed scalar curvature associated to a pseudo-hyper-K\"ahler manifold.

In particular, we know by Proposition~\ref{prop:rotating} that the result of applying the rigid c-map to a CASK manifold $(M,g,J,\nabla,\xi)$ fulfills the necessary conditions required for applying the HK/QK correspondence. In view of the results of \cite{ACDM}, the composition of these two constructions with the choice of Hamiltonian functions $f^c_Z:=f_Z -\frac12 c$ and $f^c_\h:=f_\h -\frac12 c$, where $f_Z$ and $f_\h$ are as in Proposition~\ref{prop:rotating}, will be called the \emph{supergravity c-map} in this paper.

Note that the case $c=0$ is distinguished and is called the \emph{undeformed} supergravity c-map (corresponding to \cite{FS}), while the remaining cases are collectively referred to as the \emph{deformed} supergravity c-map (corresponding to \cite{RSV}). It was shown in \cite{CHM,CDS}, under appropriate assumptions\footnote{For instance, the assumption for $c=0$ is simply that the underlying projective special K\"ahler manifold is complete, see \cite[Theorem~5]{CHM} for details. For $c>0$ the assumptions are specified in \cite[Theorems~13 and~27]{CDS}.} on the CASK manifold, that the quaternionic K\"ahler metrics $(\bar N^c,g^c_{\bar N},\mathcal Q^c)$ are complete if and only if $c\ge0$. We will therefore be assuming $c\ge 0$ henceforth, although our methods work more generally. Note that, for a fixed CASK manifold, all the manifolds $(\bar N^c,g^c_{\bar N})$ in the above family for different values of $c>0$ are locally isometric (see \cite[Proposition~10]{CDS}). For a discussion of further global properties of the deformed supergravity c-map, see \cite{MS22}.

The HK/QK correspondence, as described above, was interpreted as an instance of an even more general construction called the \emph{twist construction} in \cite{MS14}. Roughly speaking, this construction, introduced earlier by Swann, takes as input a manifold $N$ with \emph{twist data}, i.e.\ a triple $(\omega, Z, f)$ consisting of \begin{itemize}
	\item an integral closed two-form $\omega$,
	\item a vector field $Z$ generating a circle action which is Hamiltonian with respect to $\omega$,
	\item a choice of nowhere vanishing Hamiltonian function $f$,
\end{itemize} and gives as output a new manifold $\bar N$ with a circle action (and in fact, ``dual'' twist data, but this will not be important for our purposes). Furthermore, it also gives a bijective correspondence called $\mathcal H$-\emph{relatedness} between tensor fields of the same type on $N$ and $\bar N$ which are invariant under the respective circle actions. In particular, if two functions $f\in C^\infty(N)$ and $\bar f\in C^\infty(\bar N)$ are invariant under the respective circle actions and $\mathcal H$-related, then they are either both constant or both non-constant.

We refer the reader to \cite{Swann,MS14,CST21,CST22} for the details, and only summarise some of the conclusions obtained from this perspective that we will be relying on for our results.

Note that HK/QK data on a pseudo-hyper-K\"ahler manifold $N$ automatically give rise to twist data $(\omega_\h,Z,f^c_\h)$ on $N$. In fact, Macia and Swann prove the following.

\begin{thm}[{\cite[Theorem~1]{MS14}}]
	Let $(N,g_{N},I_1,I_2,I_3)$ be a pseudo-hyper-K\"ahler manifold equipped with HK/QK data $(Z,\omega_1,\omega_\h,f^c_Z,f^c_\h)$. Then the quaternionic K\"ahler manifold $(\bar N^c, g^c_{\bar N},\mathcal Q^c)$ given by the HK/QK correspondence is obtained by performing the twist construction with respect to twist data $(\omega_\h,Z,f^c_\h)$. In particular, $\mathcal Q^c$ is $\mathcal H$-related to $\Span\{I_1,I_2,I_3\}$ and $g^c_{\bar N}$ is $\mathcal H$-related to the metric $$g^c_\h:=K\left(\dfrac{1}{f^c_Z}g_N|_{{(\H Z)}^\perp}+\dfrac{f^c_\h}{(f^c_Z)^2}g_N|_{\H Z}\right),$$ where $K$ is a non-zero constant of the same sign as $f^c_Z$.
\end{thm}

Taking $K$ to have the same sign as $f^c_Z$ gives a quaternionic K\"ahler metric $g^c_{\bar N}$ that is positive-definite whenever the given pseudo-hyper-K\"ahler metric $g_N$ is positive-definite when restricted to $(\H Z)^\perp$. The reduced scalar curvature of $g^c_{\bar N}$ is then given by $\nu = -\frac{1}{8K}$. Thus, the sign of $f^c_Z$ determines the sign of the scalar curvature (they are opposite) while the choice of constant $K$ determines its magnitude. In particular, for the supergravity c-map, $f^c_Z$ is taken to be positive, so we may set $K=1$. This gives us a positive-definite supergravity c-map metric of reduced scalar curvature $-\frac{1}{8}$.

The twist construction was furthermore used in \cite{CST22} to obtain a tensor $\widetilde \Rm\in \Gamma\big((T^* N)^{\otimes 4}\big)$ on the pseudo-hyper-K\"ahler manifold $N$ to which the (lowered) Riemann curvature $\Rm_{\bar N}\in \Gamma\big((T^*\bar N^c)^{\otimes 4}\big)$ of the quaternionic K\"ahler metric $g^c_{\bar N}$ is $\mathcal H$-related. In order to state the result, we will first need to introduce some notation.

\begin{mydef}\leavevmode
	\begin{numberedlist}
		\item We define the Kulkarni--Nomizu map
		\begin{equation*}
			\begin{tikzcd}[row sep=0]
				\Gamma\big((T^*N)^{\otimes 4}\big)\ar[r] & \Gamma\big(\bigwedge^2T^*N\otimes \bigwedge^2T^*N\big)\\
				\Phi \ar[r,mapsto] & \Phi^{\owedge}
			\end{tikzcd}
		\end{equation*}
		by setting 
		\begin{equation*}
			\Phi^\owedge(A,B,C,X)\coloneqq \Phi(A,C,B,X)-\Phi(A,X,B,C)+\Phi(B,X,A,C)-\Phi(B,C,A,X)
		\end{equation*}
		for arbitrary vector fields $A,B,C,X$.
		\item We define a second map
		\begin{equation*}
			\begin{tikzcd}[row sep=0]
				\Gamma\big(\bigwedge^2T^*N\otimes\bigwedge^2T^*N\big)\ar[r] & \Gamma\big(\bigwedge^2T^*N\otimes\bigwedge^2T^*N\big)\\
				\Phi \ar[r,mapsto] & \Phi^{\obar}
			\end{tikzcd}
		\end{equation*}
		by setting
		\begin{equation*}
			\Phi^\obar (A,B,C,X)\coloneqq \Phi^\owedge(A,B,C,X)+2\Phi(A,B,C,X)+2\Phi(C,X,A,B).
		\end{equation*}
	\end{numberedlist}
\end{mydef}

For $(0,2)$-tensors $\alpha$ and $\beta$, we set $\alpha\owedge\beta\coloneqq (\alpha\otimes\beta)^\owedge$ and analogously define $\alpha\obar\beta$. Taking $\alpha$ and $\beta$ symmetric one recovers the well-known Kulkarni-Nomizu product $\alpha\owedge\beta=\beta\owedge \alpha$, which is an abstract curvature tensor, i.e.\ $(0,4)$-tensor with the symmetries of the (lowered) Riemann curvature tensor. Taking $\alpha$ and $\beta$ skew-symmetric, $\alpha\obar\beta=\beta\obar\alpha$ is precisely six times the natural projection of the tensor $\frac12 (\alpha \otimes \beta + \beta \otimes \alpha) \in \Gamma\big(\mathrm{Sym}^2\bigwedge^2T^*N \big)$ to the subspace consisting of abstract curvature tensors.

\begin{thm}[{\cite[Theorem~3.4]{CST22}}]\label{theorem:H-related_curvature}
	Let $(N,g_{N},I_1,I_2,I_3)$ be a pseudo-hyper-K\"ahler manifold equipped with HK/QK data $(Z,\omega_1,\omega_\h,f^c_Z,f^c_\h)$ and let $(\bar N^c, g^c_{\bar N},\mathcal Q^c)$ be the quaternionic K\"ahler manifold given by the HK/QK correspondence. Then the (lowered) Riemann curvature $\Rm_{\bar N}$ of the metric $g^c_{\bar N}$ is $\mathcal H$-related to the tensor
	\begin{equation}\label{eq:Rmtilde}
		\widetilde{\Rm}=\dfrac{1}{f^c_Z}\Rm_N-\dfrac{1}{f^c_Zf^c_\h}\Rm_{\mathrm{HK}}-\frac18 \Rm_{\H\mathrm{P}},
		\end{equation}
	where $\Rm_\HK$ and $\Rm_{\H\mathrm{P}}$ are defined to be
	\begin{equation*}
	\begin{aligned}
		\Rm_{\mathrm{HK}}&:=\tfrac18 \omega_\h\obar\omega_\h+\tfrac18\sum_k\omega_\h(I_k\cdot,\cdot)\owedge\omega_\h(I_k\cdot,\cdot),\\
		\Rm_{\H\mathrm{P}}&:=-g^c_\h\owedge g^c_\h-\sum_kg^c_\h(I_k\cdot,\cdot)\obar g^c_\h(I_k\cdot,\cdot).
	\end{aligned}
	\end{equation*}
\end{thm}

Note that \eqref{eq:Rmtilde} reflects a refinement of the Alekseevsky decomposition of the curvature of a quaternionic K\"ahler manifold of reduced scalar curvature $-\frac18$ arising from the HK/QK correspondence. The first two terms on the right corresponds to the part of hyper-K\"ahler type, while the last term corresponds to $-\frac18$ times the curvature of the quaternionic projective space of unit reduced scalar curvature. In particular, it follows that both $\Rm_N$ and $\Rm_\HK$ are separately $g_\h^c$-orthogonal to $\Rm_{\H\mathrm{P}}$.

As an application of Theorem~\ref{theorem:H-related_curvature}, we see that the norm of the curvature tensor $\Rm_{\bar N}$ of the metric $g^c_{\bar N}$ on the quaternionic K\"ahler side is not constant if the norm of $\widetilde\Rm$ on the pseudo-hyper-K\"ahler side is not constant. We will indeed proceed by specialising this argument to the case of the deformed supergravity c-map in the next section.

\section{Curvature formulas}

\subsection{Curvature of special K\"ahler manifolds}

We now proceed to compute the curvature of ASK manifolds.

\begin{prop}\label{prop:ASK_Riemann}
	Let $(M,g,J,\omega,\nabla)$ be an affine special K\"ahler manifold. Then the curvature $R$ of the Levi-Civita connection $D$ is $R(X,Y)=-\frac14 [\calS_X,\calS_Y]$.
\end{prop}

\begin{myproof}
	Since the Levi-Civita connection $D$ is $\nabla + \frac12 \calS$, the curvature is
	\begin{equation*}
		\begin{split}
			R(X,Y) &= [D_X,D_Y] - D_{[X,Y]}\\
			&=[\nabla_X,\nabla_Y] -\nabla_{[X,Y]} + \tfrac12[\nabla_X,\calS_Y]- \tfrac12[\nabla_Y,\calS_X]-\tfrac12\calS_{[X,Y]} + \tfrac14[\calS_X,\calS_Y] \\
			&=\tfrac12[\nabla_X,\calS_Y]- \tfrac12[\nabla_Y,\calS_X]-\tfrac12\calS_{[X,Y]} + \tfrac14[\calS_X,\calS_Y],
		\end{split}
	\end{equation*}
	where we have used that $\nabla$ is flat. On the other hand, we have $\calS_X=g^{-1}\nabla_X g$ so
	\begin{equation*}
		\begin{split}
			&[\nabla_X,\calS_Y]-[\nabla_Y,\calS_X]-\calS_{[X,Y]}\\
			&=[\nabla_X,g^{-1}\nabla_Y g]-[\nabla_Y,g^{-1}\nabla_X g]-g^{-1}\nabla_{[X,Y]}g\\
			&=\nabla_X(g^{-1})\nabla_Y g-\nabla_Y(g^{-1})\nabla_X g+g^{-1}([\nabla_X,\nabla_Y]-\nabla_{[X,Y]})g\\
			&=-g^{-1}(\nabla_Xg)g^{-1}\nabla_Y g+g^{-1}(\nabla_Yg)g^{-1}\nabla_X g=-[\calS_X,\calS_Y],
		\end{split}
	\end{equation*}
	where we have used that $\nabla$ is flat once again in the penultimate step. Putting everything together we obtain
	$$R(X,Y)=-\tfrac12 [\calS_X,\calS_Y]+\tfrac14 [\calS_X,\calS_Y]=-\tfrac14 [\calS_X,\calS_Y].$$
\end{myproof}

As a consistency check, we obtain as an immediate consequence of the above formula and Corollary~\ref{cor:S_vanish}, the well-known result that the Riemann curvature of any K\"ahler cone vanishes when applied to vector fields generating the $\C^\times$-action.

\begin{cor}\label{lemma:R(xi)=0}
	Let $(M,g,J,\nabla,\xi)$ be a conical affine special K\"ahler manifold. Then \begin{enumerate}[\normalfont(a)]
		\item $R(\xi,\cdot)\cdot=R(\cdot,\xi)\cdot=R(\cdot,\cdot)\xi=0$.
		\item $R(J\xi,\cdot)\cdot=R(\cdot,J\xi)\cdot=R(\cdot,\cdot)J\xi=0$.
	\end{enumerate}
\end{cor}

As noted earlier, the tangent space $TN=T(T^*M)$ of the total space of the cotangent bundle $N=T^*M$ of an ASK manifold $M$ can be identified, using the flat connection $\nabla$, with $\pi^*TM \oplus \pi^*T^*M$. This allows us to relate the Riemann curvature of $g_N$ to pullbacks of tensor fields defined on the base $M$.

\begin{thm}\label{theorem:curvature_N}
	Let $(M,g,J,\nabla)$ be an ASK manifold and $(N=T^*M,g_N,I_1,I_2,I_3)$ the pseudo-hyper-K\"ahler manifold given by the rigid c-map. Then the curvature tensor $\Rm_N$ of $N$ is given by \begin{align*}
		\Rm_N(A^\hor,B^\hor,C^\hor,X^\hor)&=-\tfrac14 g\big([\calS_{A^\hor},\calS_{B^\hor}]C^\hor,X^\hor\big)\\
		\Rm_N(A^\hor,B^\hor,C^\hor,X^\ver)&=\phantom{+}0\\
		\Rm_N(A^\hor,B^\hor,C^\ver,X^\ver)&=-\tfrac14 g\big([\calS_{A^\hor},\calS_{B^\hor}](C^\ver)^\sharp,(X^\ver)^\sharp\big)\\
		\Rm_N(A^\hor,B^\ver,C^\hor,X^\ver)&=\phantom{+}\tfrac12 g\big(\calS_{A^\hor}\calS_{C^\hor}(X^\ver)^\sharp,(B^\ver)^\sharp\big)\\
		&\quad+\tfrac14 g\big(\calS_{C^\hor}\calS_{A^\hor}(X^\ver)^\sharp,(B^\ver)^\sharp\big)\\
		&\quad+\tfrac14 g\big(\calS_{A^\hor}C^\hor,\calS_{(B^\ver)^\sharp}(X^\ver)^\sharp\big)\\
		&\quad-\tfrac12 \big(\nabla^2_{A^\hor,(B^\ver)^\sharp}g\big)\big(C^\hor,(X^\ver)^\sharp\big)\\
		\Rm_N(A^\hor,B^\ver,C^\ver,X^\ver)&=\phantom{+}0\\
		\Rm_N(A^\ver,B^\ver,C^\ver,X^\ver)&=-\tfrac14 g\big([\calS_{(A^\ver)^\sharp},\calS_{(B^\ver)^\sharp}](C^\ver)^\sharp,(X^\ver)^\sharp\big)
	\end{align*} where $A,B,C,X\in T_pN$, $p\in N$, and $X^\hor\in T^\hor_pN$, $X^\ver\in T^\ver_pN$ are, respectively, horizontal and vertical components. Moreover, on the right-hand side of these formulas, horizontal and vertical vectors are identified with elements of $T_{\pi (p)}M$ and $T^*_{\pi (p)}M$, respectively, and $\alpha^\sharp\in T_{\pi (p)}M$ denotes the metric dual of $\alpha\in T_{\pi (p)}^*M$.
\end{thm}

\begin{myproof}
	This general result is obtained from a straightforward but tedious computation in local coordinates $(q^i,p_j)$ on $N$ induced by local $\nabla$-affine coordinates $(q^i)$ on $M$. First one computes the Christoffel symbols of $(N,g_N)$ in terms of the Christoffel symbols of $(M,g)$ given in \eqref{eq:Christoffel}. Then one computes the curvature tensor of $(N,g_N)$ in terms of the tensor $\mathcal S$ and the curvature tensor of $(M,g)$, given in Proposition~\ref{prop:ASK_Riemann} also in terms of $\mathcal S$. One concludes by expressing the final result in a coordinate independent way using only the above intrinsic identifications and basic properties of ASK manifolds (such as the complete symmetry of $\nabla g$).
\end{myproof}

Note that the remaining components of the Riemann curvature follow from the above by symmetries of the curvature tensor and that $\nabla^2g$ coincides with $\nabla S$, where $S$ is the totally symmetric $(0,3)$-tensor which corresponds to the $(1,2)$-tensor $\mathcal{S}$.

\begin{rem}
	If the ASK manifold $M$ satisfies that $\nabla$ coincides with the Levi-Civita connection, then $\Rm_M=0$ (since $\nabla$ is flat by definition) and hence $\Rm_N=0$ by Theorem~\ref{theorem:curvature_N}.
\end{rem}

In the case where the ASK manifold $M$ is furthermore CASK, we can say something additional.

\begin{prop}\label{cor:curvature_N_HZ=0}
	Let $(M,g,J,\nabla,\xi)$ be a CASK manifold and $(N=T^*M,g_N,I_1,I_2,I_3)$ the pseudo-hyper-K\"ahler manifold given by the rigid c-map. Then the curvature tensor $\Rm_N$ is a section of the subbundle \[\textstyle \mathrm{Sym}^2\bigwedge^2 ((\H Z)^\perp)^* \oplus \bigwedge^2 ((\H Z)^\perp)^*\vee \big(((\H Z)^\perp)^*\wedge (\H Z)^*\big) \subset \mathrm{Sym}^2 \bigwedge^2T^*N,\] where we are using the isomorphism $T^*N\cong (\H Z)^* \oplus ((\H Z)^\perp)^*$ corresponding to \eqref{decomp:eq} and $\vee$ denotes the symmetric tensor product. In particular, $\Rm_N(A,B,C,X)=0$ if at least two of the vectors $A,B,C,X$ belong to $\H Z$.
\end{prop}

\begin{myproof}
	We have seen that the curvature of $N$ is completely determined by tensors on the base $M$. Under the identifications $T_p^\hor N\cong T_{\pi (p)}M$ and $T_p^\ver N \cong T^*_{\pi (p)}M$ the horizontal vector fields $Z,I_1Z$ on $N$ are identified with the vector fields $-J\xi,\xi$ on $M$, and the vertical vector fields $I_2Z,I_3Z$ with the one-forms $\xi^\flat,(-J\xi)^\flat$ (with the convention $\omega=g(J\cdot,\cdot)$). Every term in Theorem~\ref{theorem:curvature_N} can be expressed in terms of the tensor $\calS$. From Corollary~\ref{cor:S_vanish} we know that $\calS$ vanishes on $\xi$ and $J\xi$, therefore all the curvature elements are zero taking into account that $S$ and $\nabla^2g$ are totally symmetric. In fact, the total symmetry of $S=\nabla g$ was stated in Lemma~\ref{Ssymm_lem} and implies that of $\nabla S=\nabla^2g$ using that $\nabla^2_{A,B}g=\nabla^2_{B,A}g$ since $\nabla$ is flat.
\end{myproof}

\subsection{Norm of the curvature tensor}

Let us start by defining what is a locally homogeneous manifold.

\begin{mydef}
	A Riemannian manifold $(M,g)$ of dimension $n$ is called \textit{locally homogeneous} if for all $p\in M$ there exist $n$ Killing vector fields defined in a neighborhood of $p$ which are linearly independent at $p$.
\end{mydef}

Note that a function on a connected locally homogeneous Riemannian manifold which is invariant under any locally defined isometry is necessarily constant.

We will now finally show that the curvature norm of the (deformed) local c-map metric $g^c_{\bar{N}}$ associated to a CASK manifold $M$ is not constant on the manifold $\bar N^c$ unless $c=0$. We had already argued in Section~\ref{ssec:sugra-c-map} that this is equivalent to showing that the norm $\norm{\widetilde{\Rm}}_{g^c_\h}^2$ of $\widetilde{\Rm}$ with respect to $g^c_\h$, is not constant on the rigid c-map space $N=T^*M$.

In order to compute this norm, we work in a $g_N$-orthonormal frame $\{e_i,\epsilon_\mu\}$ of $TN$ that is adapted to the quaternionic distribution $\H Z$. This means that $e_i$ span the distribution $\H Z$ and $\epsilon_\mu$ span the orthogonal complement $(\H Z)^\perp$.

In terms of this frame, the norm of an abstract $(0,4)$-curvature tensor $\mathcal C$ with respect to the metric $g^c_\h$ is given by
\begin{align*}
	\norm{\mathcal C}_{g^c_\h}^2=\hat{g}^c_\h(\mathcal C,\mathcal C)&=\dfrac{(f^c_Z)^8}{(f^c_\hor)^4}\sum\mathcal C(e_i,e_j,e_k,e_\ell)^2-4\dfrac{(f^c_Z)^7}{(f^c_\hor)^3}\sum\mathcal C(\epsilon_\mu,e_j,e_k,e_\ell)^2\\
	&\quad+2\dfrac{(f^c_Z)^6}{(f^c_\hor)^2}\sum\mathcal C(\epsilon_\mu,\epsilon_\nu,e_k,e_\ell)^2+4\dfrac{(f^c_Z)^6}{(f^c_\hor)^2}\sum\mathcal C(\epsilon_\mu,e_j,\epsilon_\lambda,e_\ell)^2\\
	&\quad-4\dfrac{(f^c_Z)^5}{f^c_\hor}\sum\mathcal C(\epsilon_\mu,\epsilon_\nu,\epsilon_\lambda,e_\ell)^2+(f^c_Z)^4\sum\mathcal C(\epsilon_\mu,\epsilon_\nu,\epsilon_\lambda,\epsilon_\sigma)^2,
\end{align*}
where $\hat g^c_\h:=\big((g^c_\h)^{-1}\big)^{\otimes 4}$ denotes the metric on the bundle $(T^*N)^{\otimes 4}$ induced by $g^c_\h$.
	
Let us now specialise Theorem~\ref{theorem:H-related_curvature} to the case of the deformed supergravity c-map. Since the decomposition between the hyper-K\"ahler part and the projective quaternionic space part is orthogonal, we have $$\begin{aligned}
	\norm{\widetilde{\Rm}}_{g^c_\h}^2&=\dfrac{1}{(f^c_Z)^2}\norm{\Rm_N}_{g^c_\h}^2+\dfrac{1}{(f^c_Z)^2(f^c_\h)^2}\norm{\Rm_{\mathrm{HK}}}_{g^c_\h}^2\\
	&\quad+\dfrac{2}{(f^c_Z)^2f^c_\h}\hat{g}^c_\h(\Rm_N,\Rm_{\mathrm{HK}})+\frac{1}{64}\norm{\Rm_{\H\mathrm{P}}}_{g^c_\h}^2.
\end{aligned}$$

The final term $\frac{1}{64}\norm{\Rm_{\H\mathrm{P}}}_{g^c_\h}^2$ is a constant depending only on the dimension. Meanwhile the remaining terms can be computed to be
\begin{align*}
	\dfrac{1}{(f_Z^c)^2}\norm{\Rm_N}_{g^c_\h}^2&=\dfrac{(f_Z^c)^6}{(f_\h^c)^4}R_0^N-4\dfrac{(f_Z^c)^5}{(f_\h^c)^3}R_1^N+2\dfrac{(f_Z^c)^4}{(f_\h^c)^2}R_{2\ra}^N\\
	&\quad+4\dfrac{(f_Z^c)^4}{(f_\h^c)^2}R_{2\rb}^N-4\dfrac{(f_Z^c)^3}{f_\h^c}R_3^N+(f_Z^c)^2R_4^N,\\
	\dfrac{1}{(f_Z^c)^2(f_\h^c)^2}\norm{\Rm_{\HK}}_{g^c_\h}^2&=\dfrac{(f_Z^c)^6}{(f_\h^c)^6}R_0^{\HK}-4\dfrac{(f_Z^c)^5}{(f_\h^c)^5}R_1^{\HK}+2\dfrac{(f_Z^c)^4}{(f_\h^c)^4}R_{2\ra}^{\HK}\\
	&\quad+4\dfrac{(f_Z^c)^4}{(f_\h^c)^4}R_{2\rb}^{\HK}-4\dfrac{(f_Z^c)^3}{(f_\h^c)^3}R_3^{\HK}+\dfrac{(f_Z^c)^2}{(f_\h^c)^2}R_4^{\HK},\\
	\dfrac{1}{(f_Z^c)^2f_\h^c}\hat{g}^c_\h(\Rm_N,\Rm_{\HK})
	&=\dfrac{(f_Z^c)^6}{(f_\h^c)^5}R_0^{\mathrm C}-4\dfrac{(f_Z^c)^5}{(f_\h^c)^4}R_1^{\mathrm C}+2\dfrac{(f_Z^c)^4}{(f_\h^c)^3}R_{2\ra}^{\mathrm C}\\
	&\quad+4\dfrac{(f_Z^c)^4}{(f_\h^c)^3}R_{2\rb}^{\mathrm C}-4\dfrac{(f_Z^c)^3}{(f_\h^c)^2}R_3^{\mathrm C}+\dfrac{(f_Z^c)^2}{f_\h^c}R_4^{\mathrm C}.
\end{align*}

In the above, we have introduced the notation
$$\begin{aligned}
	R_0^N&:=\sum\Rm_N(e_i,e_j,e_k,e_\ell)^2,\\
	R_{2\ra}^N&:=\sum\Rm_N(\epsilon_\mu,\epsilon_\nu,e_k,e_\ell)^2,\\
	R_3^N&:=\sum\Rm_N(\epsilon_\mu,\epsilon_\nu,\epsilon_\lambda,e_\ell)^2,
\end{aligned}\qquad\begin{aligned}
	R_1^N&:=\sum\Rm_N(\epsilon_\mu,e_j,e_k,e_\ell)^2,\\
	R_{2\rb}^N&:=\sum\Rm_N(\epsilon_\mu,e_j,\epsilon_\lambda,e_\ell)^2,\\
	R_4^N&:=\sum\Rm_N(\epsilon_\mu,\epsilon_\nu,\epsilon_\lambda,\epsilon_\sigma)^2.
\end{aligned}$$
The terms of the form $R^\HK$ and $R^{\mathrm C}$ (where C stands for ``cross-terms'') are defined in a similar way, for example, $R^{\mathrm C}_0:=\sum\Rm_N(e_i,e_j,e_k,e_\ell)\Rm_{\HK}(e_i,e_j,e_k,e_\ell)$. In particular, all the terms $R^N$ and $R^\HK$ are non-negative functions since they are sums of squares, and by virtue of Proposition~\ref{cor:curvature_N_HZ=0}, $R^N_I=R^{\mathrm C}_I=0$ for $I=0,1,2\ra,2\rb$.

We can now express the derivative of the function $\norm{\widetilde{\Rm}}_{g^c_\h}^2$ in terms of the ``curvature functions'' $R^N, R^\HK, R^{\mathrm C}$ above.

\begin{prop}\label{crucial:prop}
	Let $(M,g,J,\nabla,\xi)$ be a CASK manifold and $(N=T^*M,g_N,I_1,I_2,I_3)$ the pseudo-hyper-K\"ahler manifold given by the rigid c-map. Let $\Xi$ be the natural lift of $\xi$ to $T^*M$ given in local $\nabla$-affine coordinates $(q^i,p_j)$ by $$\Xi:=\sum\left(q^i\dfrac{\partial}{\partial q^i}+p_i\dfrac{\partial}{\partial p_i}\right)\in\Gamma(TN).$$
	Then, the derivative of the function $\norm{\widetilde{\Rm}}_{g^c_\h}^2$ along $\Xi$ is given by \begin{equation}\label{Lie:eq}
		\calL_\Xi\norm{\widetilde{\Rm}}_{g^c_\h}^2=\dfrac{1}{(f_\hor^c)^7}\left(\sum_{k=1}^9\widetilde{\Omega}_kc^k\right),
	\end{equation} where the ($c$-independent) functions $\widetilde{\Omega}_k$ are given in terms of the functions $R^N, R^\HK, R^{\mathrm C}$ defined above by
	\begin{subequations}
		\begin{align}
			\widetilde{\Omega}_9&:=\frac{1}{128}\big(36R^N_3+R^N_4\big)\label{eq:Om9}\\
			\widetilde{\Omega}_8&:=-\frac{1}{64}\big(f_Z(260R^N_3-6R^N_4)-4R^{\mathrm{C}}_3+R^{\mathrm{C}}_4\big)\label{eq:Om8}\\
			\widetilde{\Omega}_7&:=\frac{1}{32}\big((f_Z)^2(572R^N_3+14R^N_4)+f_Z(28R^{\mathrm{C}}_3-7R^{\mathrm{C}}_4)\big)\label{eq:Om7}\\
			\widetilde{\Omega}_6&:=-\frac{1}{16}\big((f_Z)^3(788R^N_3-14R^N_4)+(f_Z)^2(-36R^{\mathrm{C}}_3+17R^{\mathrm{C}}_4)\nonumber\\
			&\quad+f_Z(-6R^{\HK}_0+20R^{\HK}_1-8R^{\HK}_{2\ra}-16R^{\HK}_{2\rb}+12R^{\HK}_3-2R^{\HK}_4)\big)\label{eq:Om6}\\
			\widetilde{\Omega}_5&:=\frac{1}{8}\big((f_Z)^4(660R^N_3)+(f_Z)^3(-36R^{\mathrm{C}}_3-15R^{\mathrm{C}}_4)\nonumber\\
			&\quad+(f_Z)^2(-30R^{\HK}_0+60R^{\HK}_1-8R^{\HK}_{2\ra}-16R^{\HK}_{2\rb}-12R^{\HK}_3+6R^{\HK}_4)\big)\label{eq:Om5}\\
			\widetilde{\Omega}_4&:=-\frac{1}{4}\big((f_Z)^5(204R^N_3+14R^N_4)+(f_Z)^4(84R^{\mathrm{C}}_3-5R^{\mathrm{C}}_4)\nonumber\\
			&\quad+(f_Z)^3(-60R^{\HK}_0+40R^{\HK}_1+16R^{\HK}_{2\ra}+32R^{\HK}_{2\rb}-24R^{\HK}_3-4R^{\HK}_4)\big)\label{eq:Om4}\\
			\widetilde{\Omega}_3&:=\frac{1}{2}\big((f_Z)^6(20R^N_3-14R^N_4)+(f_Z)^5(-12R^{\mathrm{C}}_3+19R^{\mathrm{C}}_4)\nonumber\\
			&\quad+(f_Z)^4(-60R^{\HK}_0-40R^{\HK}_1+16R^{\HK}_{2\ra}+32R^{\HK}_{2\rb}+24R^{\HK}_3-4R^{\HK}_4)\big)\label{eq:Om3}\\
			\widetilde{\Omega}_2&:=-\big((f_Z)^7(28R^N_3+6R^N_4)+(f_Z)^6(-44R^{\mathrm{C}}_3-13R^{\mathrm{C}}_4)\nonumber\\
			&\quad+(f_Z)^5(-30R^{\HK}_0-60R^{\HK}_1-8R^{\HK}_{2\ra}-16R^{\HK}_{2\rb}+12R^{\HK}_3+6R^{\HK}_4)\big)\label{eq:Om2}\\
			\widetilde{\Omega}_1&:=-2\big((f_Z)^8(8R^N_3+R^N_4)+(f_Z)^7(-20R^{\mathrm{C}}_3-3R^{\mathrm{C}}_4)\nonumber\\
			&\quad+(f_Z)^6(6R^{\HK}_0+20R^{\HK}_1+8R^{\HK}_{2\ra}+16R^{\HK}_{2\rb}+12R^{\HK}_3+2R^{\HK}_4)\big)\label{eq:Om1}
		\end{align}
	\end{subequations}
\end{prop}

\begin{myproof}
	Working in the local coordinates, we see that the vector field $\Xi$ satisfies $$\calL_\Xi g_N=2g_N,\qquad \calL_\Xi \omega_\h=2\omega_\h,\qquad\calL_\Xi f_Z=2f_Z,\qquad \calL_\Xi f_\h=2f_\h.$$
	Since we have $f^c_Z = f_Z -\frac12 c$ and $f^c_\h =f_\h -\frac12 c=-f^c_Z-c$, it follows that
	$$\calL_\Xi f^c_Z=2f^c_Z+c,\qquad \calL_\Xi f^c_\h=-2f^c_Z-c.$$
	Note in particular that $\Xi$ generates homotheties with respect to the metric $g_N$. Using the general fact that any homothety of a pseudo-Riemannian manifold is affine with respect to the Levi-Civita connection and hence preserves its curvature, we have $\calL_\Xi \Rm_N=2\Rm_N$. A straightforward computation using these observations and the Leibniz rule then yields the desired result.	
\end{myproof}

\begin{prop}\label{main:prop}
	Let $(M,g,J,\nabla,\xi)$ be a CASK manifold and $(N=T^*M,g_N,I_1,I_2,I_3)$ the pseudo-hyper-K\"ahler manifold given by the rigid c-map. The norm $\norm{\widetilde{\Rm}}_{g^c_\h}^2$ of $\widetilde{\Rm}$ with respect to $g^c_\h$ as defined in \eqref{eq:Rmtilde} is not constant on $N$ when $c>0$.
\end{prop}

\begin{myproof}
	We prove it by contradiction. Let $F^c:=\norm{\widetilde{\Rm}}_{g^c_\hor}^2\in\calC^\infty(N)$ and $\bar{F}^c:=\norm{\Rm_{\bar{N}}}_{g_{\bar{N}}^c}^2\in\calC^\infty(\bar{N})$. Suppose that $F^c$ is constant for some $c>0$. Since they are $\calH$-related, $F^c$ is constant if and only if $\bar{F}^c$ is constant. We know that for $c,c'>0$, the quaternionic K\"ahler manifolds $(\bar{N},g_{\bar{N}}^c)$ and $(\bar{N},g_{\bar{N}}^{c'})$ are locally isometric.
	
	Since there exists a (local) diffeomorphism $\varphi:\bar{N}\longrightarrow\bar{N}$ such that $\varphi^*\bar{F}^c=\bar{F}^{c'}$, it follows that $\bar{F}^c$ is constant if and only if $\bar{F}^{c'}=\varphi^*\bar{F}^c$ is constant. This implies that $\bar{F}^c$ is constant for all $c>0$. By $\calH$-relatedness, $F^c$ is also constant for all $c>0$. Then $\calL_\Xi F^c=0$ for all $c>0$ and, by \eqref{Lie:eq}, this implies that $\widetilde{\Omega}_k\equiv0$ for $k=1,\ldots,9$.
	
	By \eqref{eq:Om9}, $\widetilde{\Omega}_9\equiv0$ implies that $$36R^N_3+R^N_4\equiv0,$$ but both functions are non-negative, so this means that $R^N_3\equiv R^N_4\equiv0$. Recall that $R^N_4$ is a sum of squares, so each of the individual terms must vanish separately, i.e. $\Rm_N(\epsilon_\mu,\epsilon_\nu,\epsilon_\lambda,\epsilon_\sigma)\equiv0$. This shows that $\Rm_N\equiv0$, which implies $R^{\mathrm{C}}_3\equiv0$ and $R^{\mathrm{C}}_4\equiv0$.
	
	Now, by \eqref{eq:Om1}, $\widetilde{\Omega}_1\equiv0$ implies that $$6R^{\HK}_0+20R^{\HK}_1+8R^{\HK}_{2\ra}+16R^{\HK}_{2\rb}+12R^{\HK}_3+2R^{\HK}_4\equiv0,$$ but, as before, all these functions are non-negative, so all of them vanish identically. Thus, we find that $\Rm_{\HK}\equiv0$, but this is a contradiction, since for a rigid c-map space we have $\Rm_{\HK}(Z,I_1Z,Z,I_1Z)=g_N(Z,Z)^2>0$. Hence we can conclude that $\calL_\Xi\norm{\widetilde{\Rm}}_{g^c_\h}^2\not\equiv0$ and therefore $\norm{\widetilde{\Rm}}_{g^c_\h}^2$ is not a constant function.
\end{myproof}

As a consequence, we obtain our main theorem.

\begin{thm}\label{main:thm}
	Let $(M,g,J,\nabla,\xi)$ be a CASK manifold and $(\bar N^c,g^c_{\bar N},\mathcal Q^c)$ the quaternionic K\"ahler manifold given by the deformed supergravity c-map. Then, for any $c>0$, $(\bar N^c,g^c_{\bar N})$ is not locally homogeneous.
\end{thm}

By the results of \cite{CST21,CRT,MS22}, given a CASK manifold $(M,g,J,\omega,\nabla,\xi)$ of real dimension $2n$ with automorphism group $\Aut(M)$, by which we mean that $\Aut(M)$ is the group of isometries of $(M,g)$ preserving the full CASK data, the associated supergravity c-map spaces $(\bar N^c,g^c_{\bar N},\mathcal Q^c)$ are isometrically acted on by the group $\Aut(M)\ltimes\Heis_{2n+1}$, provided that the underlying projective special K\"ahler manifold $(\bar M, g_{\bar M})$ is simply connected\footnote{The assumption of simply connectedness can be dropped if $M$ is a CASK domain.}. In particular, when $\Aut(M)$ acts transitively on the underlying projective special K\"ahler manifold $\bar M$, we obtain an action of $\Aut(M)\ltimes \Heis_{2n+1}$ that is transitive on the level sets of the norm of the quaternionic moment map associated to the circle action on $\bar N$. Thus, as a corollary of our main theorem, we have the following generalisation of a result in \cite{CST22} concerning the supergravity c-map space associated to a flat CASK manifold.

\begin{cor}\label{main:cor}
	Let $(M,g,J,\omega,\nabla,\xi)$ be a CASK manifold fibering over a simply connected projective special K\"ahler manifold $(\bar M, g_{\bar M})$. Assume that its automorphism group $\Aut(M)$ acts transitively on $\bar M$. Then the associated supergravity c-map space $(\bar N^c,g^c_{\bar N},\mathcal Q^c)$ is a complete quaternionic K\"ahler manifold of cohomogeneity exactly one when $c>0$.
\end{cor}

\begin{myproof}
	Since the Riemannian manifold $(\bar M, g_{\bar M})$ is complete, the corresponding undeformed c-map space $(\bar N,g_{\bar N})=(\bar N^0,g^0_{\bar N})$ is complete in virtue of \cite[Theorem~10]{CHM}. As a first step, we will show that $(\bar N,g_{\bar N})$ is not only complete but is in fact homogeneous.
	
	The group of isometries $\Aut(\bar M)\subset \Isom(\bar M)$ induced by $\Aut(M)$ extends canonically to a group of isometries of $(\bar N,g_{\bar N})$. This is stated in \cite[Proposition~26]{CDJL} for CASK domains but holds in general as a consequence of \cite[Lemma~4]{CHM}. It can be also seen as a special case ($c=0$) of the results of \cite{CST21,CRT,MS22} mentioned above. The group $\Aut(\bar M)$ acts transitively on the base of the fiber bundle $\bar N \rightarrow\bar M$ mapping fibers to fibers. In addition, there is a fiber-preserving isometric action of the solvable Iwasawa subgroup $G_{2n+2}$ of $\mathrm{SU}(1,n+1)$ on $\bar N|_{\bar U}$ \cite[Theorem~5]{CHM} for every domain $U\subset M$, which is isomorphic to a CASK domain, where $\bar U$ denotes the image of $U$ under the projection $M\rightarrow \bar M$. (Recall that every CASK manifold is locally isomorphic to a CASK domain.) Note that $\dim G_{2n+2}=2n+2$, where $\dim \bar M = 2n-2$.
	This solvable group action on $\bar N|_{\bar U}$ is simply transitive on each fiber. In particular, for every such $\bar U$ there is a Lie algebra $\mathfrak{g}_{\bar U}\cong \mathfrak{g}_{2n+2}= \mathrm{Lie}(G_{2n+2})$ of Killing fields of $\bar N|_{\bar U}$ transitive on each fiber. Moreover, $\mathfrak{g}_{\bar U}$ can be identified with the space of parallel sections over $\bar U$ of a flat symplectic vector bundle over $\bar M$ (with Lie algebras as fibers), compare \cite[Theorem~9]{CHM}. Since $M$ is simply connected the above vector bundle has a global parallel frame. Thus we obtain a globally defined Lie algebra of Killing fields $\mathfrak{g}\cong\mathfrak{g}_{2n+2}$ of $\bar N$, which restricts to $\mathfrak{g}_{\bar U}$ on the domain $\bar N|_{\bar U}\subset \bar N$. Since $\bar N$ is complete, there is a corresponding Lie group $G$ acting on $\bar N$, which together with $\Aut(\bar M)$ generates a transitive group of isometries of $\bar N$.
	
	Now that we know that $(\bar N,g_{\bar N})$ is a \emph{homogeneous} quaternionic K\"ahler manifold of negative scalar curvature, we can apply the following arguments to show that it belongs to the class of Alekseevsky spaces. First of all, in virtue of the resolution of the Alekseevsky conjecture about the structure of homogeneous Einstein manifolds of negative scalar curvature by B\"ohm and Lafuente \cite{BL}, we know that $(\bar N,g_{\bar N})$ admits a simply transitive solvable Lie group of isometries. Then, by a result of Lauret \cite{Lauret}, it is a standard Einstein solvmanifold in the sense of Heber \cite{H}. Finally, by \cite[Theorem~B]{H} such a manifold admits a simply transitive completely solvable group of isometries. Quaternionic K\"ahler manifolds with that property were classified by Alekseevsky \cite{A,Co}.
	
	We claim that the one-loop deformation of any c-map space which is an Alekseevsky space is complete if the deformation parameter $c$ is positive (for $c=0$ it holds by homogeneity). First we note that all of the Alekseevsky spaces with exception of the quaternionic hyperbolic spaces and the Hermitian symmetric spaces of non-compact type dual to complex Grassmannians of $2$-planes can be represented as q-map spaces \cite{dW}, a special class of c-map spaces. By \cite[Theorem~27]{CDS} the one-loop deformation of a complete q-map space is complete if $c>0$. In particular, the one-loop deformed Alekseevsky q-map spaces with $c>0$ are complete. Furthermore, the Hermitian symmetric Alekseevsky spaces were shown to have regular boundary behavior, implying the completeness of their one-loop deformation for $c>0$ \cite[Example~14 and Theorem~13]{CDS}.
	
	Finally, we are left with the quaternionic hyperbolic spaces $\HH^n$. We claim that these cannot be represented as c-map spaces and hence cannot occur in our setting. This can be seen by looking at totally geodesic K\"ahler submanifolds compatible with the quaternionic structure. Thanks to \cite{AlekM} we know that the maximal possible dimension of a K\"ahler submanifold compatible with the quaternionic structure of a quaternionic K\"ahler manifold of dimension $4n$ is $2n$. In the case of $\HH^n$ the only totally geodesic K\"ahler submanifolds of (real) dimension $2n$ compatible with the quaternionic structure are the complex hyperbolic subspaces $\CH^n$ (up to isometries of the ambient space). On the other hand, any c-map space of dimension $4n$ has a totally geodesic K\"ahler submanifold compatible with the quaternionic structure of the form $\CH^1\times\bar M$, where $\bar M$ is the underlying projective special K\"ahler manifold of dimension $2n-2$. In fact, the submanifold $\CH^1 \times \bar M\subset \bar N$ is obtained as the fixed point set of the isometric involution expressed in standard fiber coordinates $(\rho, \tilde\phi, \tilde\zeta_i,\zeta^i)$, $i=1,\ldots,2n$, \cite{CHM} by $(\rho, \tilde\phi, \tilde\zeta_i,\zeta^i)\mapsto (\rho, \tilde\phi, -\tilde\zeta_i,-\zeta^i)$. Since $\CH^n$ is irreducible, we see that $\HH^n$ is not a c-map space if $n>1$. The case $n=1$ is also excluded, since the c-map space associated with a projective special K\"ahler manifold reduced to a point is $\CH^2$ (belonging to the Hermitian symmetric series) and not $\HH^1$. This finishes the proof of the completeness of $(\bar N^c,g^c_{\bar N})$ for $c>0$.
	
	Now the corollary follows from the fact that $(\bar N^c,g^c_{\bar N})$ has a group of isometries acting with cohomogeneity one but no such group acting with cohomogeneity $0$.
\end{myproof}

\bibliographystyle{alpha}
\bibliography{References}

\end{document}